\documentclass[12pt]{article}

\usepackage{IEEEtrantools}
\usepackage{amssymb}
\usepackage{amsmath,amsfonts,amsthm}
\usepackage{mathrsfs}
\usepackage{url}
\usepackage{tikz}

\newcommand{\bbbt}{\mathbb{T}}
\newcommand{\scrt}{\mathscr{T}}

\newcommand{\be}{\begin{equation}}
\newcommand{\ee}{\end{equation}}
\newcommand{\bea}{\begin{eqnarray}}
\newcommand{\eea}{\end{eqnarray}}
\newcommand{\bean}{\begin{eqnarray*}}
\newcommand{\eean}{\end{eqnarray*}}
\newcommand{\brray}{\begin{array}}
\newcommand{\erray}{\end{array}}
\newcommand{\biearray}{\begin{IEEEarray}{rCl}}
\newcommand{\eiearray}{\end{IEEEarray}}

\newcommand{\newsection}[1]{\setcounter{equation}{0}
\setcounter{dfn}{0}
\section{#1}}

\newcommand{\newsubsection}[1]{
\subsection{#1}}

\newtheorem{dfn}{Definition}[section]
\newtheorem{thm}[dfn]{Theorem}
\newtheorem{lmma}[dfn]{Lemma}
\newtheorem{ppsn}[dfn]{Proposition}
\newtheorem{crlre}[dfn]{Corollary}
\newtheorem{xmpl}[dfn]{Example}
\newtheorem{rmrk}[dfn]{Remark}

\newcommand{\bdfn}{\begin{dfn}\rm}
\newcommand{\bthm}{\begin{thm}}
\newcommand{\blmma}{\begin{lmma}}
\newcommand{\bppsn}{\begin{ppsn}}
\newcommand{\bcrlre}{\begin{crlre}}
\newcommand{\bxmpl}{\begin{xmpl}}
\newcommand{\brmrk}{\begin{rmrk}\rm}

\newcommand{\edfn}{\end{dfn}}
\newcommand{\ethm}{\end{thm}}
\newcommand{\elmma}{\end{lmma}}
\newcommand{\eppsn}{\end{ppsn}}
\newcommand{\ecrlre}{\end{crlre}}
\newcommand{\exmpl}{\end{xmpl}}
\newcommand{\ermrk}{\end{rmrk}}

\newcommand{\bbc}{\mathbb{C}}
\newcommand{\bbz}{\mathbb{Z}}
\newcommand{\bbn}{\mathbb{N}}
\newcommand{\bbr}{\mathbb{R}}

\newcommand{\clh}{\mathcal{H}}

\newcommand{\clk}{\mathcal{K}}
\newcommand{\cll}{\mathcal{L}}

\newcommand{\cls}{\mathcal{S}}

\newcommand{\prf}{\noindent{\it Proof\/}: }

\def \qed { \mbox{}\hfill
$\Box$\vspace{1ex}}

\begin{document}


 \author{\sc  {Bipul Saurabh}}
 \title{ $q$-invariance of quantum quaternion spheres}
 \maketitle


\begin{abstract} 
The $C^*$-algebra of continuous functions on  the quantum quaternion sphere $H_q^{2n}$ 
can  be identified with the quotient algebra $C(SP_q(2n)/SP_q(2n-2))$. 
In commutative case i.e. for $q=1$, the topological 
space $SP(2n)/SP(2n-2)$ is homeomorphic to the odd dimensional sphere $S^{4n-1}$. In this paper, we prove the noncommutative analogue of this result.
Using homogeneous $C^*$-extension theory, we prove that the $C^*$-algebra $C(H_q^{2n})$ is isomorphic to the $C^*$-algebra  $C(S_q^{4n-1})$. This
  further implies that for different values of $q \in [0,1)$, the $C^*$-algebras underlying the noncommutative space $H_q^{2n}$ are isomorphic.
\end{abstract}

{\bf AMS Subject Classification No.:} {\large 58}B{\large 34}, {\large
46}L{\large 80}, {\large
  19}K{\large 33} 
\newsection{Introduction}
Quantization of Lie groups and their homogeneous spaces have played an important role in 
linking the theory of  compact quantum group with  noncommutative geometry.  
Many authors (see \cite{VakSou-1990ab},  \cite{PodVai-1999aa}, \cite{ChaPal-2008aa}, 
 \cite{PalSun-2010aa}) have studied different aspects of the theory of quantum homogeneous spaces. However, in these 
papers, main examples have been the quotient spaces of the compact quantum group $SU_q(n)$.  
 Neshveyev \& Tuset
(\cite{NesTus-2012ab}) studied  quantum homogeneous spaces in a more general set up and 
gave a complete classification of the irreducible representations
of the $C^*$-algebra $C(G_q/H_q)$ where $G_q$ is the $q$-deformation of a simply
connected semisimple compact Lie group and $H_q$ is the $q$-deformation of a closed
Poisson-Lie subgroup $H$ of $G$. Moreover, Neshveyev \& Tuset 
(\cite{NesTus-2012ab}) proved that $C(G_q/H_q)$ is $KK$-equivalent to the classical counterpart
$C(G/H)$. Quantum symplectic group $SP_q(2n)$ and its homogeneous space $C(SP_q(2n)/SP_q(2n-2))$ have been studied by 
the author in \cite{Sau-2015aa} and 
 $K$-groups of quotient space $C(SP_q(2n)/SP_q(2n-2))$ with explicit generators were obtained.

The $C^*$-algebra $C(H_q^{2n})$ of  continuous functions on the quantum quaternion sphere is defined as the universal $C^*$-algebra 
given by a finite set of generators and relations (see \cite{Sau-2015aa}). In \cite{Sau-2015aa},
the isomorphism between the quotient algebra 
$C(SP_q(2n)/SP_q(2n-2))$ and $C(H_q^{2n})$ has been established. Now several questions arise about this noncommutative space $H_q^{2n}$.
\begin{enumerate}
\item
Is $H_q^{2n}$ topologically same as $S_q^{4n-1}$, i.e.\  are the $C^*$-algebras
$C(H_q^{2n})$ and $C(S_q^{4n-1})$ isomorphic?
\item
Are the $C^*$-algebras $C(H_q^{2n})$ isomorphic for different values of $q$?
\item
Does the quantum quaternion sphere admit a good
spectral triple equivariant under the $SP_q(2n)$-group action?
\end{enumerate}
 We attempt the first two questions in this paper.  In commutative case i.e. for $q=1$, 
the quotient space $SP(2n)/SP(2n-2)$ can be realized as the  quaternion sphere $H^{2n}$.
It can be easily verified that the quaternion sphere $H^{2n}$
is homeomorphic to the odd dimensional sphere $S^{4n-1}$. One can now expect the quotient algebra $C(SP_q(2n)/SP_q(2n-2))$
or equivalently the $C^*$-algebra $C(H_q^{2n})$ to be isomorphic to the 
$C^*$-algebra underlying the odd dimensional quantum sphere $S_q^{4n-1}$. In this paper, using homogeneous $C^*$-extension theory, 
we show that this is indeed the case.

The remarkable work done by L. G. Brown, R. G. Douglas and P. A. Fillmore (\cite{BroDouFil-1977aa}) on extensions 
of commutative $C^*$-algebras by compact operators has 
led many authors to extend this theory further in order to provide a tool for analysing the structure of $C^*$-algebras. For a nuclear, separable 
$C^*$-algebra $A$ and a separable $C^*$-algebra $B$, G. G. Kasparov (\cite{Kas-1979aa}) 
constructed the group $Ext(A,B)$ consisting ``stable equivalence classes'' of $C^*$-algebra extensions of the form 
\[
 0 \rightarrow B\otimes \clk \rightarrow E \rightarrow A \rightarrow 0
\]
Here $E$ will be called the middle $C^*$-algebra.
One of the important features of this construction is that the group $Ext(A,B)$ 
coincides with the group $KK^1(A,B)$. Another important aspect is that  it does not demand much.
It does not require the extensions to be unital or essential. But at the same time, it does not provide much 
informations about the  middle $C^*$-algebras. 
Since elements of the group $Ext(A,B)$ are stable equivalence classes and not unitary equivalence classes of extensions, 
two elements in the same class 
may have nonisomorphic middle $C^*$-algebras. For a nuclear $C^*$-algebra $A$ and a
finite dimensional compact metric space $Y$ (i.e. a closed subset of $\cls^n$ for some $n \in \bbn$), 
M. Pimsner, S. Popa and D. Voiculescu~(\cite{PimPopVoi-1979aa}) constructed another group $Ext_{PPV}(Y,A)$ 
consisting of unitary equivalence classes of unital homogeneous 
extensions of $A$ by $C(Y) \otimes \clk$. For $y_0 \in Y$, the subgroup $Ext_{PPV}(Y,y_0,A)$ consists of those elements of $Ext_{PPV}(Y,A)$ 
that split at $y_0$. For a commutative $C^*$-algebra $A$, the group  $Ext_{PPV}(Y,A)$ was computed by Schochet in \cite{Sch-1980aa}. Further 
Rosenberg \& Schochet~(\cite{RosCho-1981ab}) showed that $Ext_{PPV}(Y,A^+)=Ext(A,C(Y))$ and $Ext_{PPV}(Y^+,+,A^+)=Ext(A,C(Y))$ where $Y$ is a 
finite dimensional locally compact Hausdorff space, $+$ is the point at infinity and $A^+$ is the $C^*$-algebra obtained by adjoining unity to $A$.

To show that the $C^*$-algebra $C(H_q^{2n})$ is isomorphic to $C(S_q^{4n-1})$, 
we first exhibit an isomorphism between the group $Ext_{PPV}(Y, y_0,A)$ and the group $Ext_{PPV}(Y, y_0,\Sigma^2A)$ under certain assumptions 
on the topological space $Y$
where $\Sigma^2A$ is the quantum double suspension of $A$ and $y_0 \in Y$. Using this, we describe all elements of  the group
$Ext_{PPV}(\bbbt,C(S_0^{2\ell+1}))$ explicitly. We then prove that all nonisomprphic 
middle $C^*$- algebras that occur in all the extensions  of the group  $Ext_{PPV}(\bbbt,C(S_0^{2\ell+1}))$ have different  $K$-groups. 
Then using representation theory of $C(H_q^{2n})$, 
we show that the  following extension 
\[
 0 \rightarrow C(\bbbt)\otimes \clk \rightarrow C(H_q^{2n}) \rightarrow C(S_0^{4n-3}) \rightarrow 0
 \]
 is unital and homogeneous. Now by comparing the $K$-groups, we prove that the above extension is unitarily equivalent to either the following  
 extension
 \[
 0 \rightarrow C(\bbbt)\otimes \clk \rightarrow C(S_0^{4n-1}) \rightarrow C(S_0^{4n-3}) \rightarrow 0
 \]
 or  its inverse in the group $Ext_{PPV}(\bbbt,C(S_0^{2\ell+1}))$. 
 This proves that the $C^*$-algebras $C(H_q^{2n})$ and $C(S_0^{4n-1})$ are isomorphic. For $q=0$, it follows 
 directly from the defining relations. In \cite{PalSun-2010aa}, it was proved that for different values of $q \in [0,1)$
 the $C^*$-algebras $C(S_q^{4n-1})$ are isomorphic. As a consequence,  the 
 $C^*$-algebras $C(H_q^{2n})$ and $C(S_q^{4n-1})$ are isomorphic for all $q \in [0,1)$.
 Also, this shows that  the $C^*$-algebras $C(H_q^{2n})$ are isomorphic for different values of $q$ which establishes  $q$-invariance 
 of the quantum quaternion spheres.

We now set up some notations.
    The standard
bases of the Hilbert spaces $L_2(\bbn)$ and $L_2(\bbz)$ will be denoted by
$\left\{e_n: n\in \bbn \right\}$ and $\left\{e_n: n\in \bbz \right\}$ respectively.  
We  denote the left shift operator on $L_2(\bbn)$  and $L_2(\bbz)$ by the same notation $S$. For $m< 0$, $(S^*)^m$ 
denotes the operator $S^{-m}$. Let $p_i$ denote the rank one projection sending 
$e_i$ to $e_i$ and $p$ denote the operator $p_0$. We write $\cll(\clh)$ and $\clk(\clh)$ for the sets of 
all bounded linear operators on $\clh$ and compact operators on $\clh$  respectively. We denote by  $\clk$  the $C^*$-algebra 
of compact operators. For a $C^*$-algebra $A$, $\Sigma^2A$ and $M(A)$ are 
used to denote the quantum double suspension of $A$ and multiplier algebra of $A$ respectively. 
The map  $\pi$ will  denote the canonical homomorphism from $M(A)$ to $Q(A):=M(A)/A$ and for $a \in M(A)$, 
$[a]$ stands for the image of $a$ under the map $\pi$. For a locally compact Hausdorff space $Y$, 
we write $Y^+$ to denote one point compactification of $Y$. For a $C^*$-algebra $A$, $A^+$ denotes 
the $C^*$-algebra obtained by adjoining unity to $A$. Both the symbols $\cls^n$ and  $\bbbt^n$ will denote 
the $n$-dimensional sphere. 
Unless otherwise stated,  $q$ will denote a
real number in the interval $(0,1)$.

\newsection { $C^*$-algebra extensions }
We first  recall some notions related to  the $C^*$-extension theory. Let $A$ be a unital separable  nuclear  $C^*$-algebra. 
Let  $B$ be a stable $C^*$- algebra. An extension of 
$A$ by $B$ is a
short exact sequence $0 \rightarrow B \stackrel{i}{\rightarrow} E \stackrel{j}{\rightarrow} A \rightarrow 0$. In such case, there exists a unique 
homomorphism $\sigma: E \rightarrow M(B)$ such that  $\sigma(i(b))=b$ for all $b \in B$. 
We can now define the Busby invariant for the extension $0 \rightarrow B \stackrel{i}{\rightarrow} E \stackrel{j}{\rightarrow} A \rightarrow 0$
as the homomorphism 
$\tau : A \rightarrow M(B)/B$ given by $ \tau(a)=\pi \circ \sigma (e)$ where $e$ is a preimage of $a$ and $\pi$ is the canonical map 
$M(B) \rightarrow M(B)/B$. It is easy to see that $\tau$ is well defined.  An
extension $\tau$ is called an essential extension  if $\tau$ is injective or equivalently image of $B$ is an essential ideal of $E$. 
We call an extension unital if it is a unital homomorphism or equivalently $E$ is a unital $C^*$-algebra.
An extension $\tau$ is called a trivial (or split) extension if there exists a homomorphism $\lambda:A \rightarrow M(B)$ such that 
$\tau=\pi \circ \lambda$. Two extensions $\tau_1$ and $\tau_2$ are said to be   weakly unitarily  equivalent if there exists a unitary $u$ in $Q(B)$  
such that $u\tau_1(a)u^*=\tau_2(a)$ for all $a \in A$. 
They are said to be  unitarily  equivalent if there exists a unitary $U$ in $M(B)$ 
 such 
that  $\pi(U)\tau_1(a)\pi(U^*)=\tau_2(a)$ for all $a \in A$.  We denote unitarily equivalence relation by $\sim_u$. 
Let $Ext_{\sim_u}(A,B)$ denote the set of unitary equivalence classes of extensions of $A$ by $B$. One can put a
binary operation ''+`` on  $Ext_{\sim_u}(A,B)$ as follows. 
Since $M(B)$ is a stable $C^*$-algebra, we can get two isometries $\nu_1$ and $\nu_2$ in $M(B)$ such that $\nu_1\nu_1^*+\nu_2\nu_2^*=1$. 
Let $\tau_1$ and $\tau_2$ be two elements in $Ext_{\sim_u}(A,B)$. Define $\tau_1+\tau_2: A \rightarrow Q(B)$ by 
\begin{IEEEeqnarray}{rCl} \label{sec1-+}
 (\tau_1+\tau_2)(a):=\pi(\nu_1)\tau_1(a)\pi(\nu_1^*)+\pi(\nu_2)\tau_2(a)\pi(\nu_2^*).
\end{IEEEeqnarray}
This makes $Ext_{\sim_u}(A,B)$ a commutative semigroup. Moreover,  the set of trivial extensions 
forms a subsemigroup of $Ext_{\sim_u}(A,B)$. We denote the quotient of $Ext_{\sim_u}(A,B)$ with the set of trivial extensions
 by $Ext(A,B)$. 
 For  a separable nuclear $C^*$-algebra $A$, the set 
$Ext(A,B)$ under the operation $+$ is a group (see \cite{Bla-1998aa}).
Two extensions $\tau_1$ and $\tau_2$ represent  the same element in $Ext(A,B)$ if 
 there exists two trivial extensions $\phi_1$ and $\phi_2$ such that $\tau_1+\phi_1 \sim_u \tau_2 +\phi_2$. 
 One can show that for a stable $C^*$-algebra $B$, $Ext(A,B)=Ext(A,B\otimes \clk)$. Now for an
arbitrary $C^*$-algebra $B$,  define $Ext(A,B):=Ext(A,B\otimes \clk)$. 
We denote an equivalent class in the group $Ext(A,B)$ of an extension $\tau$ by $[\tau]_s$.
For $B=\bbc$, we denote the group $Ext(A,\bbc)$ 
by $Ext(A)$. Note that in this case, two unital essential extensions $\tau_1$ and $\tau_2$ are in the same equivalence
class (i.e. $[\tau_1]_s=[\tau_2]_s$) 
if and only if they are unitarily equivalent.

Suppose that  $Y$ is a  finite dimensional compact metric space i.e. a closed subset of $\cls^n$ for some $n \in \bbn$. 
Let $M(Y)$, $Q(Y)$ and $Q$ be 
the $C^*$-algebras $M(C(Y)\otimes \clk)$, $M(C(Y)\otimes \clk)/C(Y)\otimes \clk$ and $\cll(\clh)/\clk(\clh)$ (Calkin algebra) respectively.
It 
is easy to show that $M(Y)$ is the set of all continuous functions from $Y$ to $\cll(\clh)$  where continuity 
is with respect to $*$-strong operator topology on $\cll(\clh)$. 
We call an extension $\tau$ of $A$ by $C(Y) \otimes \clk$ homogeneous if for all $y \in Y$,
the map $ev_y \circ \tau : A \rightarrow Q$ is injective 
 where $ev_y:Q(Y) \rightarrow Q$ is the evaluation map at $y$. 
Let $Ext_{PPV}(Y,A)$ be the set of unitary equivalence classes of unital homogeneous extensions of $A$ by $C(Y)\otimes \clk$. For a 
nuclear $C^*$-algebra $A$, Pimsner, Popa and Voiculescu (\cite{PimPopVoi-1979aa}) showed that 
$Ext_{PPV}(Y,A)$ is a group with the additive operation defined in \ref{sec1-+}. We denote the equivalence class in the group 
$Ext_{PPV}(Y,A)$ of an extension $\tau$ by $[\tau]_u$.  For $y_0 \in Y$, define the set 
\[
 Ext_{PPV}(Y,y_0,A)=\left\{[\tau]_u \in Ext_{PPV}(Y,A): ev_{y_0} \circ \tau \mbox{ is trivial} \right\}.
\]
The set $Ext_{PPV}(Y,y_0,A)$ is a subgroup of $Ext_{PPV}(Y,A)$.

\newsubsection{The groups $Ext_{PPV}(Y, A)$ and $Ext_{PPV}(Y, \Sigma^2A)$}
Here we will show that for a separable nuclear $C^*$-algebra $A$ and a finite dimensional compact metric 
space $Y$ such that $K$-groups of $C(Y)$ are free groups with finite generators, 
the groups   $Ext_{PPV}(Y, A)$ and $Ext_{PPV}(Y,\Sigma^2A)$ are isomorphic.  Let us recall some definitions. 
We say that two elements $a$ and $b$ in $Q(A)$ are  unitarily equivalent  if there exists a unitary 
$U \in M(A)$ such that 
$[U]a[U^*]=b$.  They are weakly unitarily equivalent if there exists unitary $u \in Q(A)$ such that $uau^*=b$. 
We call an element $a$ in a $C^*$-algebra $B$ norm-full if 
it is not contained in any proper closed ideal in $B$. Suppose that $A$ and $B$ are separable $C^*$-algebras. An extension
$\tau : A \rightarrow Q(B\otimes \clk)$ is said to be norm-full if for every nonzero element $a \in A$, $\tau(a)$ is norm full 
element of $Q(B\otimes \clk)$.
\bdfn Let $B$ be a separable stable $C^*$-algebra. Then $B$ is said to have the corona factorization
property if every norm-full projection in $M(B)$ is Murray-von Neumann equivalent to unit element of $M(B)$.
\edfn
It is easy to see that a $C^*$-algebra $A$ with 
corona factorization
property, any norm-full projection in $Q(B)$ is Murray-von Neumann equivalent to $1$ of $Q(B)$. Further, one can show that for
a finite dimensional compact metric 
space $Y$, $C(Y)\otimes \clk$ has corona factorization property (see \cite{PimPopVoi-1979aa}).
\bppsn \label{injectivekk1}
Let $A$ be a unital separable nuclear $C^*$-algebra which satisfies the Universal Coefficient
Theorem. Suppose that  $Y$ is  a finite dimensional compact metric 
space. Then the map  
\begin{IEEEeqnarray}{rCl}
 i: Ext_{PPV}(Y, A) & \longrightarrow & KK^1(A,C(Y)) \nonumber \\
   \quad  [\tau]_u & \mapsto & [\tau]_s \nonumber 
\end{IEEEeqnarray}
is an injective homomorphism. 
\eppsn
\prf
Since unitarily equivalence implies stable equivalence, the map $i$ is well defined.   Any unital homogeneous extension
  is a purely large extension and  hence a norm-full extension (see page 19, \cite{EllKuc-2001aa}). 
  Therefore  from Theorem 2.4 in \cite{Hua-2004aa}, it follows  that $i$ is injective.
\qed

From now on, 
without  loss of generality, we will assume that 
the Hilbert space $\clh$ is $L_2(\bbn)$. 
Let $\tau$ be a unital homogeneous extension of $A$ by $C(Y)\otimes \clk(\clh)$.
Define $\tilde{\tau} : A \rightarrow Q(C(Y)\otimes \clk(\clh) \otimes \clk(\clh))$ by
: $\tilde{\tau}(a)=[\tau_a\otimes p]$ where $[\tau_a]=\tau(a)$. By universal property of quantum double suspension
(see proposition 2.2, \cite{HonSzy-2008aa}), we have a homomorphism 
\begin{IEEEeqnarray}{rCl} \label{sigmatau}
\Sigma^2 \tau : \Sigma^2A &\rightarrow & Q(C(Y)\otimes \clk(\clh) \otimes \clk(\clh)) 
\end{IEEEeqnarray}
such that 
$\Sigma^2 \tau(a \otimes p)=\tilde{\tau}(a)=[\tau_a \otimes p]$ and $\Sigma^2 \tau(1 \otimes S)= [1 \otimes 1\otimes S]$. 
Clearly $\Sigma^2 \tau$ is a unital extension. Since $\tau$ is homogeneous, the map $ev_y \circ \Sigma^2 \tau$ 
is injective on the $C^*$- algebra $A\otimes p$ for all $y \in Y$. Making use of the fact that
$(1\otimes p)A \otimes \clk(1\otimes p)=A\otimes p$, 
one can prove that the map $ev_y \circ \Sigma^2 \tau$ is 
injective on $A \otimes \clk$. Since $A\otimes \clk$ is an essential ideal of $\Sigma^2A$, we conclude that the map $ev_y \circ \Sigma^2 \tau$ 
is injective on $\Sigma^2A$ and hence 
 $\Sigma^2 \tau$ is a homogeneous extension.  Moreover, if $\tau_1$ and $\tau_2$ are unitarily equivalent by a unitary $U \in M(C(Y)\otimes \clk(\clh))$ 
 then so are $\Sigma^2 \tau_1$ and $\Sigma^2 \tau_2$ by the unitary $U\otimes 1 \in M(C(Y)\otimes \clk(\clh)\otimes \clk(\clh))$. 
 This gives a well defined map 
 \begin{IEEEeqnarray}{rCl}
 \beta: Ext_{PPV}(Y, A) & \longrightarrow & Ext_{PPV}(Y,\Sigma^2A) \nonumber \\
   \quad  [\tau]_u & \mapsto & [\Sigma^2 \tau]_u \nonumber 
 \end{IEEEeqnarray}
 \bppsn \label{injectivesigma}
 The map $\beta: Ext_{PPV}(Y,  A)  \longrightarrow Ext_{PPV}(Y,  \Sigma^2A)$ given above is injective group homomorphism.
 \eppsn
 \prf
 It follows from straightforward calculations.
 \qed\\
 To get surjectivity of the map $\beta$, we need to put certain assumptions on the topological space $Y$.
 \bppsn \label{visometry}
 Let $Y$ be a finite dimensional compact metric space. Assume that 
 $K_0(C(Y))$ and $K_1(C(Y))$ are free groups with finite number of generators. Let   $V \in Q(C(Y)\otimes \clk(\clh)\otimes \clk(\clh))$ be a
 isometry such that $VV^*$ and $1-VV^*$ both
 are norm full  projections. Then $V$  
 is weakly unitarily equivalent to $[1\otimes1\otimes S]$.
 \eppsn
 \prf
 We assume that $V$ is not weakly 
 unitarily equivalent to $[1\otimes1\otimes S]$.
 Since $C(Y)\otimes \clk$ has corona factorization property, it follows that 
 $VV^*$ and $1-VV^*$ both are
 Murray-von Neumann equivalent to $[1]$ of $Q(C(Y)\otimes \clk(\clh) \otimes \clk(\clh))$. 
 Also, one can easily verify that $[1\otimes 1\otimes p]$ and $[1-1\otimes 1\otimes p]=[1 \otimes 1 \otimes (1-p)]$ are 
 Murray-von Neumann equivalent to $[1]$ of $Q(C(Y)\otimes \clk(\clh) \otimes \clk(\clh))$. This implies that $VV^*$ is weakly unitarily 
 equivalent to $1-[1\otimes 1\otimes p]$. So, without loss of generality, we can assume that $V$ has final projection $1-[1\otimes 1\otimes p]$. 
 Take a split unital homogeneous extension $\tau$ of $C(\bbbt)$ by $C(Y)\otimes \clk(\clh)$. Let $\Sigma_V^2\tau$ be  a  unital homogeneous 
 extension  of 
 $\Sigma^2C(\bbbt)$ by  $C(Y)\otimes \clk(\clh\otimes \clk(\clh)$ 
 given by $\Sigma_V^2 \tau(a \otimes p)=[\tau_a \otimes p]$ and $\Sigma_V^2 \tau(1 \otimes S)= V$ where $[\tau_a]=\tau(a)$. 
  From Corollary 3.8 (\cite{Hua-2004aa}) and  the fact that 
 $V$ is not weakly 
 unitarily equivalent to $[1\otimes1\otimes S]$, it follows  
 that $[\Sigma_V^2\tau]_u$ is not in the image of the map $\beta$.
 Let $n[\Sigma_V^2 \tau]_u= \beta([\phi]_u)$ for some $n \in \bbz-\{0\}$ and for some 
 unital homogeneous extension $\phi$ of $C(\bbbt)$ by $C(Y)\otimes \clk(\clh)$. It is easy to see that $\phi$ must be split and in that case 
 $n[\Sigma_V^2 \tau]$ is the class of split extensions. By proposition \ref{injectivekk1} and 
 the fact that $KK^1(\Sigma^2C(\bbbt),C(Y))$ is free group, we get that 
 $\Sigma_V^2 \tau$ is a  split extension. This contradicts the fact that $[\Sigma_V^2\tau]_u$ is not in the image of the map $\beta$. So, for any 
 $n \in \bbz-\{0\}$,  $n[\Sigma_V^2 \tau]_u$ is not in the image of the map $\beta$.  This shows that image of 
 $Ext_{PPV}(Y,\Sigma^2C(\bbbt))$ in the group $KK^1(\Sigma^{2}C(\bbbt),C(Y))$ has one more free generator than the group $Ext_{PPV}(Y,C(\bbbt))$ in 
 $KK^1(C(\bbbt),C(Y))\equiv KK^1(\Sigma^{2}C(\bbbt),C(Y))$.
 Since for all $n \in \bbn$,  $KK^1(\Sigma^{2n}C(\bbbt),C(Y))\equiv K_0(C(Y))\oplus K_1(C(Y))$ are  free groups with finite generators, 
 iterating 
 above process will lead to a contradiction. This proves that $V$ is  weakly 
 unitarily equivalent to $[1\otimes1\otimes S]$.
 \qed 
 \brmrk
 Here we should point out that above proposition may hold for any finite dimensional compact metric space $Y$. 
 But since we could not find it in literature, we prove the proposition under certain assumptions on $Y$.
 \ermrk
 
 \bcrlre \label{vunitarilyequi}
 Let $Y$ and $V$ be as in the above proposition. Then $V$  
 is  unitarily equivalent to $[1\otimes1\otimes S]$.
  \ecrlre
  \prf 
  Consider the unital extension $\Sigma_V^2\tau$ constructed in proposition \ref{visometry}. From Corollary 3.8 (\cite{Hua-2004aa}), it follows that 
  $\Sigma_V^2\tau$ is unitarily equivalent to $\Sigma^2\tau$ defined in equation \ref{sigmatau} with $A=C(\bbbt)$. Hence $V$  
 is  unitarily equivalent to $[1\otimes1\otimes S]$.
  \qed 
 
The following lemma establishes the isomorphism between the groups $Ext_{PPV}(Y,  A)$ and $Ext_{PPV}(Y,  \Sigma^2A)$ under certain assumptions 
on the space  $Y$.
 
 \blmma \label{QDS}
 Let $Y$ be a finite dimensional compact metric space. Assume that the groups 
 $K_0(C(Y))$ and $K_1(C(Y))$ are free groups with finite number of generators.
 Then the map $\beta: Ext_{PPV}(Y,  A)  \longrightarrow Ext_{PPV}(Y,  \Sigma^2A)$ given above is an isomorphism.
 \elmma
 
 \prf
 We only need to show that $\beta$ is surjective thanks to proposition \ref{injectivesigma}. 
 Let $\phi$ be a unital homogeneous extension of $\Sigma^2A$ by $C(Y)\otimes \clk(\clh) \otimes \clk(\clh)$. Let $\phi(1\otimes S)=V$. Since 
 $\phi$ is a unital homogeneous extension and hence a norm full extension, it follows that 
 $VV^*$ and $1-VV^*$ are norm full projections. Therefore by Corollary \ref{vunitarilyequi}, there exists a unitary
 $U \in M(C(Y)\otimes \clk(\clh) \otimes \clk(\clh))$ such that $[U]V[U^*]=[1\otimes 1\otimes S]$. 
 So without loss of generality, we can assume that $\phi$ maps  $1\otimes S$ to $[1\otimes1 \otimes  S]$. 
 This implies that $\phi(1\otimes p)=[1\otimes 1\otimes p]$. 
But then $\phi(A\otimes p) \subset (1\otimes p)\phi(A\otimes p)(1\otimes p) \subset  Q(C(Y)\otimes \clk(\clh)) \otimes p$ 
 which induces a map $\tau :A \rightarrow Q(C(Y)\otimes \clk(\clh))$  by omitting the projection $p$. 
 Therefore we get a unital homogeneous  extension of $A$ 
 such that $\beta([\tau]_u)=[\phi]_u$. Hence $\beta$ is surjective. 
 \qed 
  \bcrlre \label{QDSsplit}
  For $y_0 \in Y$, 
  the map 
  \[\beta_{|Ext_{PPV}(Y, y_0,  A)}: Ext_{PPV}(Y, y_0,  A)  \longrightarrow Ext_{PPV}(Y, y_0,  \Sigma^2A)
   \]
 is an isomorphism.
  \ecrlre
  \prf
  It is easy to check that if $ev_{y_0} \circ \tau$ is split then so is 
 $ev_{y_0} \circ \Sigma^2 \tau$ and vice versa. Now the claim will follow by Lemma \ref{QDS}.
 \qed

\newsection{ Elements of $Ext_{PPV}(\bbbt,C(S_0^{2\ell+1}))$ }

In the present  section, we will write down   all elements of the groups $Ext(C(S_0^{2\ell+1}))$ and 
$Ext_{PPV}(\bbbt,C(S_0^{2\ell+1}))$ in terms of their Busby invariants.
Define the 
$*$-homomorphisms $\varphi_m$
 as follows:
\begin{IEEEeqnarray}{rCl}
 \varphi_m: C(S_0^{2\ell+1}) &\rightarrow & 
 Q\Big(\clk \big(\underbrace{L_2(\bbn)\otimes \cdots \otimes L_2(\bbn)}_{\ell+1 \mbox{ copies } }\big)\Big) \nonumber \\
 S^*\otimes 1 \otimes \cdots \otimes 1 &\mapsto & S^*\otimes 1 \otimes \cdots \otimes 1 \nonumber \\
 p\otimes S^*\otimes 1 \otimes \cdots \otimes 1 &\mapsto & p\otimes S^*\otimes 1 \otimes \cdots \otimes 1 \nonumber \\
 \cdots \nonumber \\
 p\otimes p \otimes \cdots \otimes p \otimes S^*\otimes 1 &\mapsto & p\otimes p \otimes \cdots \otimes S^* \otimes 1  \nonumber \\
 p\otimes p \otimes \cdots \otimes t &\mapsto &   p\otimes p \otimes \cdots \otimes p \otimes (S^*)^m \nonumber
\end{IEEEeqnarray}
It is easy to verify that $\varphi_m$'s are essential unital extensions of $C(S_0^{2\ell+1})$ by compact operators. Hence 
$[\varphi_m]_s \in Ext(C(S_0^{2\ell+1}))$. We shall show that each element in the group $Ext(C(S_0^{2\ell+1}))$  is of the form  
$[\varphi_m]_s$ for some $m \in \bbz$.
Let $\clh_0$ be the Hilbert space $\underbrace{L_2(\bbn)\otimes \cdots \otimes L_2(\bbn)}_{\ell \mbox{ copies } }\otimes L_2(\bbz)$.  For 
$m \in \bbz$, let 
$\vartheta_m$ be the representation of $C(S_0^{2\ell+1})$ given by 
\begin{IEEEeqnarray}{rCl}
 \vartheta_m: C(S_0^{2\ell+1}) &\rightarrow & 
 \cll(\clh_0) \nonumber \\
 S^*\otimes 1 \otimes \cdots \otimes 1 &\mapsto & S^*\otimes 1 \otimes \cdots \otimes 1 \nonumber \\
 p\otimes S^*\otimes 1 \otimes \cdots \otimes 1 &\mapsto & p\otimes S^*\otimes 1 \otimes \cdots \otimes 1 \nonumber \\
 \cdots \nonumber \\
 p\otimes p \otimes \cdots \otimes p \otimes S^*\otimes 1 &\mapsto & p\otimes p \otimes \cdots \otimes S^* \otimes 1  \nonumber \\
 p\otimes p \otimes \cdots \otimes t &\mapsto &   p\otimes p \otimes \cdots \otimes p \otimes (S^*)^m\nonumber
\end{IEEEeqnarray}
Let $P$ be the self adjoint projection in $\cll(\clh_0)$ on 
the subspace spanned by the basis elements  
$\left\{e_{n_1} \otimes \cdots \otimes e_{n_{\ell+1}}: n_i \in \bbn \mbox { for all } i \in \left\{1,2, \cdots ,\ell+1\right\}\right\}$. 
One can check that 
$\mathcal{F}_m:=\Big(C(S_0^{2\ell+1}), \clh_0, 2P-1 \Big)$ 
with the underlying representation $\vartheta_m$ is a Fredholm module. By Proposition 17.6.5 in (\cite{Bla-1998aa}, page 157), 
 the group $Ext(C(S_0^{2\ell+1}))$ is isomorphic to the group $K^1(C(S_0^{2\ell+1}))$. Under this identification, one can easily show that 
 the equivalence class of the Fredholm module $\mathcal{F}_m$ 
 corresponds to the equivalence class  $[\varphi_m]_s$. 

\bppsn For $\ell \in \bbn$, one has
\[
 Ext(C(S_0^{2\ell+1}))=\left\{[\varphi_m]_s: m \in \bbz \right\}.
\]
\eppsn
\prf  To prove the claim, we will use the index pairing between the groups $K_1(C(S_0^{2\ell+1}))$ 
and $K^1(C(S_0^{2\ell+1}))$ given by Kasparov product (see \cite{Bla-1998aa}). The group $K_1(C(S_0^{2\ell+1}))$ is 
generated  by the unitary 
$u:=\underbrace{p\otimes \cdots \otimes p}_{\ell \mbox{ copies} } \otimes t + 1 - \underbrace{p\otimes \cdots \otimes p}_{\ell \mbox{ copies} } \otimes 1$.
 For $m \in \bbz$, let  $R_m : P\clh_0 \rightarrow P\clh_0$ be the operator $P\vartheta_m(u)P$.
Hence we get
\[
 \langle u, \mathcal{F}_m \rangle = \mbox{Index}(R_m)=m.
\]
This completes the proof.
\qed \\
To describe all elements of $Ext_{PPV}(\bbbt, C(S_0^{2\ell+1}))$, we define the 
$*$-homomorphisms $\phi_m$
 as follows:
\begin{IEEEeqnarray}{rCl}
 \phi_m: C(S_0^{2\ell+1}) &\rightarrow & 
 Q\Big(\clk \big(\underbrace{L_2(\bbn)\otimes \cdots \otimes L_2(\bbn)}_{\ell+1 \mbox{ copies } }\big)\otimes C(\bbbt)\Big) \nonumber \\
 S^*\otimes 1 \otimes \cdots \otimes 1 &\mapsto & S^*\otimes 1 \otimes \cdots \otimes 1 \nonumber \\
 p\otimes S^*\otimes 1 \otimes \cdots \otimes 1 &\mapsto & p\otimes S^*\otimes 1 \otimes \cdots \otimes 1 \nonumber \\
 \cdots \nonumber \\
 p\otimes p \otimes \cdots \otimes p \otimes S^*\otimes 1 &\mapsto & p\otimes p \otimes \cdots \otimes S^* \otimes 1 \otimes 1 \nonumber \\
 p\otimes p \otimes \cdots \otimes t &\mapsto &   p\otimes p \otimes \cdots \otimes p \otimes (S^*)^m \otimes 1 \nonumber
\end{IEEEeqnarray}
It is easy to verify that $\phi_m$'s are essential unital extensions. Since last component  is $1$, these extensions are also homogeneous. 
Let $A_m$ be the $C^*$-subalgebra of $C(S_0^{2\ell+3})$ generated by  the operators 
\begin{IEEEeqnarray}{lCl}
 S^*\otimes 1 \otimes \cdots \otimes 1 \otimes 1, \nonumber \\
 p\otimes S^*\otimes 1 \otimes \cdots \otimes 1 \otimes 1 \nonumber \\
 \qquad \qquad \cdots \nonumber \\
 p\otimes p \otimes \cdots \otimes S^* \otimes 1 \otimes 1, \nonumber \\
  p\otimes p \otimes \cdots \otimes p \otimes (S^*)^m \otimes 1 \nonumber 
\end{IEEEeqnarray}
and
$\clk \big(\underbrace{L_2(\bbn)\otimes \cdots \otimes L_2(\bbn)}_{\ell+1 \mbox{ copies } }\big)\otimes C(\bbbt)$. 
Then for each $m \in \bbz$, we have the following exact sequence
\[
 0\longrightarrow \clk \big(\underbrace{L_2(\bbn)\otimes \cdots \otimes L_2(\bbn)}_{\ell+1 \mbox{ copies } }\big)\otimes C(\bbbt) 
 \longrightarrow A_m \longrightarrow C(S_0^{2\ell+1}) \longrightarrow 0
\]
with the Busby invariant $\phi_m$. By using the six term sequence, one can show that 
\begin{IEEEeqnarray}{rCl} \label{K-groups}
 K_0(A_m)=\bbz \oplus \bbz/m\bbz, \qquad \quad K_1(A_m)=\bbz.
\end{IEEEeqnarray}
\blmma \label{split} For  $\ell \in \bbn$ and $t_0 \in \bbbt$, one has
 \[
  Ext_{PPV}\big(\bbbt, t_0, C(S_0^{2\ell+1})\big)=\left\{ 0\right\}, \qquad Ext_{PPV}\big(\bbbt,C(S_0^{2\ell+1})\big)=\bbz.
 \]
 \elmma
 \prf
 It follows from Theorem 1.5 in \cite{RosCho-1981ab} that 
 \[
  Ext_{PPV}\big(\bbbt, t_0, C(\bbbt)\big)= Ext_{PPV}\big(\bbr^+, t_0, C_0(\bbr)^+\big)=Ext\big(C_0(\bbr),C_0(\bbr)\big)=\left\{0\right\}.
 \]
The $C^*$-algebra $C(S_0^{2\ell+1})$ can be obtained by applying quantum double suspension on $C(\bbbt)$ repeatedly (see \cite{HonSzy-2002ab}). 
Therefore from Corollary \ref{QDSsplit}, 
 we have
\[
 Ext_{PPV}\big(\bbbt, t_0, C(S_0^{2\ell+1})\big)=Ext_{PPV}\big(\bbbt, t_0, C(\bbbt)\big)=\left\{0\right\}.
\] 
Further from Theorem 1.4 in \cite{RosCho-1981ab}, we get 
\[
 Ext_{PPV}\big(\bbbt,C(\bbbt)\big)=Ext_{PPV}\big(\bbbt,C_0(\bbr)^+\big)=Ext\big(C_0(\bbr),C(\bbbt)\big)=\bbz.
\]
Hence by applying Lemma \ref{QDS}, we get the claim.
 \qed \\
The following lemma says that each element of the group $Ext_{PPV}(\bbbt,C(S_0^{2\ell+1}))$ is of the form $[\phi_m]_u$ for some $m \in \bbz$.
\blmma \label{ext computation}
For $\ell \in \bbn$, one has
\[
 Ext_{PPV}(\bbbt,C(S_0^{2\ell+1}))=\left\{[\phi_m]_u: m \in \bbz \right\}.
\]
\elmma
\prf Fix $t_0 \in \bbbt$. 
Define a homomorphism $\Psi$ as follows:
\begin{IEEEeqnarray}{rCl}
 \Psi: Ext_{PPV}\big(\bbbt,C(S_0^{2\ell+1})\big) &\longrightarrow & Ext \big(C(S_0^{2\ell+1})\big) \nonumber \\
 \quad [\tau]_u & \longmapsto & [ev_{t_0} \circ \tau]_s \nonumber
\end{IEEEeqnarray} 
Clearly $\ker \Psi=Ext_{PPV}\big(\bbbt,t_0,C(S_0^{2\ell+1})\big)=\left\{0\right\}$.   Therefore $\Psi$ is an injective group homomorphism.
Since for all $m \in \bbz$,  $ev_{t_0}\circ \phi_m=\varphi_m$, it follows that the homomorphism $\Psi$ is surjective. This proves the claim.
\qed

\newsection{Quantum quaternion sphere}
In this section, we first recall the definition and representation theory  of the $C^*$-algebra $C(H_q^{2n})$ of continuous functions 
on the quantum quaternion sphere. 
Then we prove our main result 
that the $C^*$-algebra
$C(H_q^{2n})$ is isomorphic to the $C^*$-algebra  $C(S_q^{4n-1})$.
\bdfn  \label{dfn-relation}
 The $C^{*}$-algebra $C(H_{q}^{2n})$ of continuous functions on the quantum quaternion sphere is defined  as the universal  
 $C^{*}$-algebra generated by elements $z_{1}$, $z_{2}$, ....$z_{2n}$ satisfying the following relations:
\begin{align}
z_{i}z_{j} &=  qz_{j}z_{i} & \text{for } & i > j,  i+j \neq 2n+1 \label{c1}\\
z_{i}z_{i^{'}}  &= q^2z_{i^{'}}z_{i} -(1-q^{2})\sum_{k > i} q^{i-k}z_{k}z_{k^{'}} & \text{for } &  i > n \label{c2}\\
z_{i}^{*}z_{i^{'}} &= q^{2}z_{i^{'}}z_{i}^{*} \label{c3}\\
z_{i}^{*}z_{j} &= qz_{j}z_{i}^{*} & \text{for }  & i+j>2n+1,  i \neq j \label{c4}\\
z_{i}^{*}z_{j} &=qz_{j}z_{i}^{*} + (1-q^{2})\epsilon_{i}\epsilon_{j}q^{\rho_{i}+\rho_{j}}z_{i^{'}}z_{j^{'}}^{*} 
& \text{for }  & i+j<2n+1,  i \neq j \label{c5}\\
z_{i}^{*}z_{i} &=z_{i}z_{i}^{*} + (1-q^{2})\sum_{k>i}z_{k}z_{k}^{*} & \text{for } & i >n \label{c6}\\
z_{i}^{*}z_{i} &=z_{i}z_{i}^{*} +(1-q^{2})q^{2\rho_{i}}z_{i^{'}}z_{i^{'}}^{*}+ (1-q^{2})\sum_{k>i}z_{k}z_{k}^{*} 
   & \text{for }  & i \leq n \label{c7}\\
\sum_{i=1}^{2n}z_{i}z_{i}^{*} &= 1 \label{c8}
\end{align} 
\edfn
In \cite{Sau-2015aa}, we showed that the $C^*$-algebra $C(H_{q}^{2n})$ is isomorphic to the quotient algebra $C(SP_{q}(2n)/SP_{q}(2n-2))$
that can also be described as the 
$C^*$-subalgebra of $C(SP_q(2n))$ generated by $\left\{u_m^1,u_m^{2n} : m \in \{1,2,\cdots 2n\}\right\}$ i.e. elements of first and last row 
of fundamental  matrix of $C(SP_q(2n))$.
 Here we briefly describe 
all irreducible representations of $C(H_{q}^{2n})$. For a detailed treatment on this, we refer the reader to  \cite{Sau-2015aa}.
Let $N$ be the number operator given by $N:e_{n} \mapsto n e_{n}$ 
and $S$ be the  shift operator given by $S:e_{n} \mapsto e_{n-1}$ on $L_{2}(\bbn)$.
For $i=1,2,\cdots,n-1$, let $\pi_{s_{i}}$ denote the following representation of $C(SP_q(2n))$,
\[
\pi_{s_{i}}(u_l^k)=\begin{cases}
              \sqrt{1-q^{2N+2}}S & \mbox{ if } (k,l)=(i,i) \mbox{ or } (2n-i,2n-i),\cr
              S^*\sqrt{1-q^{2N+2}} & \mbox{ if } (k,l)=(i+1,i+1) \mbox{ or } (2n-i+1,2n-i+1),\cr
							-q^{N+1} & \mbox{ if } (k,l)=(i,i+1),\cr
							q^N & \mbox{ if } (k,l)=(i+1,i),\cr
							q^{N+1} & \mbox{ if } (k,l)=(2n-i,2n-i+1),\cr
							-q^N & \mbox{ if } (k,l)=(2n-i+1,2n-i),\cr
							\delta_{kl} & \mbox{ otherwise }. \cr
							\end{cases}
\]
For $i=n$,
\[
\pi_{s_n}(u_l^k)=\begin{cases}
              \sqrt{1-q^{4N+4}}S & \mbox{ if } (k,l)=(n,n),\cr
              S^*\sqrt{1-q^{4N+4}} & \mbox{ if } (k,l)=(n+1,n+1),\cr
							-q^{2N+2} & \mbox{ if } (k,l)=(n,n+1),\cr
							q^{2N} & \mbox{ if } (k,l)=(n+1,n),\cr
							\delta_{kl} & \mbox{ otherwise }. \cr
              \end{cases}
\]

Each $\pi_{s_{i}}$ is an  irreducible representation and is called an elementary representation of $C(SP_{q}(2n))$.
For any two representations $\varphi$ and $\psi$ of $C(SP_{q}(2n))$ define, $\varphi * \psi := (\varphi \otimes \psi)\circ \Delta$ where $\Delta$ 
is the co-multiplication map of $C(SP_{q}(2n))$.
Let W be the Weyl group of $sp_{2n}$ and $\vartheta \in W$ such that  $s_{i_{1}}s_{i_{2}}...s_{i_{k}}$ is a reduced expression for $\vartheta$. 
Then $\pi_{\vartheta}= \pi_{s_{i_{1}}}*\pi_{s_{i_{2}}}*\cdots *\pi_{s_{i_{k}}}$ is an irreducible representation which is independent 
of the reduced expression. Now for $t=(t_{1},t_{2},\cdots ,t_{n}) \in \bbbt^{n}$, define the map  $\tau_t: C(SP_{q}(2n) \longrightarrow \bbc $ by
\[
\tau_{t}(u_j^i)=\begin{cases}
              \overline{t_{i}}\delta_{ij} & \mbox{ if } i \leq n,\cr
							t_{2n+1-i}\delta_{ij} & \mbox{ if } i > n,\cr
							\end{cases}
\]
Then $\tau_{t}$ is a $*$-algebra homomorphism. For $t \in \bbbt^{n}, \vartheta \in W$,  let $\pi_{t,\vartheta} = \tau_{t}*\pi_{\vartheta}$. Define 
the representation 
$\eta_{t,\vartheta}$ of  $C(H_{q}^{2n})$ as  $\pi_{t,\vartheta}$ restricted to $C(H_{q}^{2n})$. 
Denote by $\omega_{k}$ the following word of Weyl group of $sp_{2n}$,
\[
\omega_k =\begin{cases}
             I & \mbox{ if } k =1, \cr
             s_1s_2\cdots s_{k-1} & \mbox{ if } 2 \leq k \leq n,\cr
             s_1s_2\cdots s_{n-1}s_{n}s_{n-1}\cdots s_{2n-k+1} & \mbox{ if } n < k \leq 2n.\cr
						\end{cases}
\] 

For $k=1$, define $\eta_{t,I}:C(H_{q}^{2n})\rightarrow \bbc \mbox{ such that } \eta_{t,I}(z_j) = t\delta_{1j}$.
The set $\left\{\eta_{t,I}: t \in T\right\}$ gives all one dimensional irreducible representations of $C(H_{q}^{2n})$. 
\bthm \label{repn} (\cite{Sau-2015aa})
The set $\left\{\eta_{t,\omega_k}: 1\leq k \leq 2n, t \in \bbbt \right\}$ gives a complete list of irreducible representations of $C(H_{q}^{2n})$.
\ethm  
Define $\eta_{\omega_k}: C(H_{q}^{2n})\rightarrow C(\bbbt) \otimes \scrt^{\otimes k-1}$ 
such that $\eta_{\omega_k}(a)(t) = \eta_{t,\omega_k}(a)$ for all $a \in  C(H_{q}^{2n})$. 
Let $C_{1}^{2n}= C(\bbbt)$ and for $ 2 \leq k \leq 2n$,  $C_{k}^{2n} = \eta_{\omega_{k}}(C(H_{q}^{2n}))$. 

\bcrlre \label{cr5}
The set $\left\{\eta_{t,\omega_l}: 1\leq l \leq k, t \in \bbbt \right\}$ gives a complete list of irreducible representations of $C_k^{2n}$. 
\ecrlre
\bcrlre
$\eta_{\omega_{k}}$ is a faithful representation of $C_k^{2n}$.
\ecrlre

By Corollary \ref{cr5}, one can find all primitive ideals i.e. kernels of  irreducible representations of $C_k^{2n}$.  Define
$y_l^{k}:=\eta_{\omega_k}(z_l)$ for $1\leq l \leq k$. Let 
$I_l^k$ be the ideal of $C_k^{2n}$ generated by $\left\{ y_l^k,y_{l+1}^k, \cdots ,y_k^k \right\}$. For  $t \in \bbbt$, let $C_t(\bbbt)$ be the 
set of all continuous functions on $\bbbt$ vanishing at the point $t$. Then 
\begin{IEEEeqnarray}{rCl} \label{ideal structure}
  \left\{C_t(\bbbt)\otimes \clk(L_2(\bbn))^{\otimes(k-1)} \right\}_{t \in \bbbt} 
 \subset I_k^k \subset I_{k-1}^k \subset \cdots \subset I_1^k=C_k^{2n}
\end{IEEEeqnarray}
is a complete list of primitive ideals of $C_k^{2n}$. In Lemma 5.1  in \cite{Sau-2015aa}, we established the following exact sequence
\[
      0 \longrightarrow C(\bbbt) \otimes \clk(L_2(\bbn))^{\otimes(k)} \longrightarrow 
      C_{k+1}^{2n} \stackrel{\sigma_{k+1}}{\longrightarrow} C_{k}^{2n} \longrightarrow 0
\] 
where $\sigma_{k+1}$ is the restriction of $(1^{\otimes (k)}\otimes \sigma)$ to $C_{k+1}^{2n}$, the map $\sigma : \scrt \rightarrow \bbc$ 
is the homomorphism such that $\sigma(S)=1$ and $\scrt$ is the toeplitz algebra.
The following lemma says that this exact sequence is a unital homogeneous extension of $C_{k}^{2n}$ by 
$C(\bbbt)\otimes \clk$.
\blmma \label{lmma-homogeneous}
For $1 \leq k \leq 2n$, the following exact sequence 
\[
      0 \longrightarrow C(\bbbt) \otimes \clk(L_2(\bbn))^{\otimes(k)} \longrightarrow 
      C_{k+1}^{2n} \stackrel{\sigma_{k+1}}{\longrightarrow} C_{k}^{2n} \longrightarrow 0
\] 
is a unital homogeneous extension of $C_{k}^{2n}$ by 
$C(\bbbt)\otimes \clk$.
\elmma
\prf Since  $C_{k+1}^{2n}$ is unital, the given extension is unital.
Let $\tau: C_k^{2n} \rightarrow Q(\bbbt)$ be the Busby invariant corresponding to this extension. For $t_0 \in \bbbt$, 
let $\tau_{t_0} : C_k^{2n} \rightarrow Q$ be
the map $ev_{t_0} \circ \tau$ where $ev_{t_0}:Q(\bbbt) \rightarrow Q$ is the evaluation map at $t_0$. Assume that 
$J_{t_0}=\ker(\tau_{t_0})$. To show that 
the given short exact sequence is a homogeneous extension, we need to prove that $J_{t_0}=\left\{0\right\}$ for all $t_0 \in \bbbt$.\\
\textbf{Case 1: $n < k < 2n$} \\
We have
\begin{IEEEeqnarray}{rCl} 
\tau_{t_0}(y_k^k)&=&\tau_{t_0}(t \otimes \underbrace{q^N \otimes \cdots \otimes q^N}_{(n-1)\mbox{ copies }} \otimes q^{2N} \otimes
\underbrace{q^N \otimes \cdots \otimes q^N}_{(k-n-1) \mbox{ copies }})\nonumber \\
&=&t_0[ \underbrace{q^N \otimes \cdots \otimes q^N}_{(n-1)\mbox{ copies }} \otimes q^{2N} \otimes
\underbrace{q^N \otimes \cdots \otimes q^N}_{(k-n-1) \mbox{ copies }} \otimes \sqrt{1-q^{2N}}S^*] \label{eqn 1} \\
&\neq & 0. \nonumber  
\end{IEEEeqnarray}
This shows $y_k^k \notin J_{t_0}$.
Since $J_{t_0}$ is intersection of all primitive ideals that contains $J_{t_0}$, we conclude that $J_{t_0}$ is equal 
to $C_F(\bbbt)\otimes \clk$ for some closed 
subset $F$ of $\bbbt$ where $C_F(\bbbt)$ is set of all continuous functions on $\bbbt$ vanishing on $F$. From equation \ref{eqn 1}, we get 
\begin{IEEEeqnarray}{rCl}
 \tau_{t_0}((y_k^k)(y_k^k)^*)&=& [\underbrace{q^{2N} \otimes \cdots \otimes q^{2N}}_{(n-1)\mbox{ copies }} \otimes q^{4N} \otimes
\underbrace{q^{2N} \otimes \cdots \otimes q^{2N}}_{(k-n-1) \mbox{ copies }} \otimes (1-q^{2N})] \nonumber \\
&=& [\underbrace{q^{2N} \otimes \cdots \otimes q^{2N}}_{(n-1)\mbox{ copies }} \otimes q^{4N} \otimes
\underbrace{q^{2N} \otimes \cdots \otimes q^{2N}}_{(k-n-1) \mbox{ copies }} \otimes 1]. \nonumber 
\end{IEEEeqnarray}
Therefore 
\begin{IEEEeqnarray}{rCl}
 \tau_{t_0}(1\otimes \underbrace{ p \otimes \cdots p}_{(k-1) \mbox{ copies }})&=&
 [\underbrace{p\otimes \cdots \otimes p}_{ (k-1) \mbox{ copies }} \otimes 1]. \nonumber 
\end{IEEEeqnarray}
Hence 
\begin{IEEEeqnarray}{rCl}
 \tau_{t_0}(t\otimes \underbrace{ p \otimes \cdots \otimes p}_{(k-1) \mbox{ copies }}) &=&
 t_0[\underbrace{p\otimes \cdots \otimes p}_{ (k-1) \mbox{ copies }} \otimes \sqrt{1-q^{2N}}S^*] \nonumber \\
 &=&t_0[\underbrace{p\otimes \cdots \otimes p}_{ (k-1) \mbox{ copies }} \otimes S^*]. \nonumber
 \end{IEEEeqnarray}
Consider the function $\chi: C(\bbbt) \rightarrow Q$ such that $\chi(t)=[S^*]$. Since $[S^*]$ is unitary in $Q$ with spectrum equal to $\bbbt$, 
it follows that the map $\chi$ is injective. 
This shows  that for any nonzero function $f$ on $\bbbt$, 
$\tau_{t_0}(f(t)\otimes \underbrace{ p \otimes \cdots \otimes p}_{(k-1) \mbox{ copies }}) \neq 0$ 
which further implies that $F=\bbbt$ and $J_{t_0}=\left\{0\right\}$. \\
\textbf{Case 2: $1 \leq k \leq n$} \\
For $k=n$, 
\begin{IEEEeqnarray}{rCl}
 \tau_{t_0}(y_n^n) &=& t_0[ \underbrace{q^N \otimes \cdots \otimes q^N}_{(n-1)\mbox{ copies }} \otimes \sqrt{1-q^{4N}}S^*]. \nonumber 
\end{IEEEeqnarray}
For $1\leq k < n$,
\begin{IEEEeqnarray}{rCl}
 \tau_{t_0}(y_k^k) &=& t_0[ \underbrace{q^N \otimes \cdots \otimes q^N}_{(k-1)\mbox{ copies }} \otimes \sqrt{1-q^{2N}}S^*]. \nonumber 
\end{IEEEeqnarray}
Similar calculations as done in the case $1$ shows that $J_{t_0}= \left\{0\right\}$. This establishes the claim.
\qed\\
We now state the main result of this paper. 
\bthm \label{isomorphism}
For all $n\in \bbn, n \geq 2$ and $1\leq k \leq 2n$, the $C^*$-algebra $C_k^{2n}$ is isomorphic to  the $C^*$-algebra $C(S_0^{2k-1})$ of continuous 
functions on
odd dimensional quantum sphere. In particular, $C(H_q^{2n})$ is 
isomorphic to $C(S_0^{4n-1})$ or equivalently to $C(S_q^{4n-1})$.
\ethm
\prf 
 Fix $n$. To prove the theorem, we use induction on $k$.  For $k=1$, $C_1^{2n}=C(\bbbt)$. So the claim is true for $k=1$. Assume that the claim is 
 true for $k$ i.e. $C_k^{2n}$ is isomorphic to 
 $C(S_0^{2k-1})$. From Lemma \ref{lmma-homogeneous}, it follows that  following short exact sequence 
\begin{IEEEeqnarray}{rCl} \label{exact}
 0 \longrightarrow C(\bbbt)\otimes \clk \longrightarrow C_{k+1}^{2n} \longrightarrow C_k^{2n} \longrightarrow 0
\end{IEEEeqnarray}
is a unital homogeneous extension. Therefore it  can be viewed  as an element of the 
 group $Ext_{PPV}(\bbbt,C(S_0^{2k-1}))$. This implies that it is unitarily equivalent to $\phi_m$ or equivalently  to the following 
 exact sequence 
\begin{IEEEeqnarray}{rCl} 
 0 \longrightarrow C(\bbbt)\otimes \clk \longrightarrow A_m \longrightarrow C(S_0^{2k-1}) \longrightarrow 0 \nonumber
\end{IEEEeqnarray}
for some $m \in \bbz$. 
 From Theorem 5.3  in \cite{Sau-2015aa} and equation $(\ref{K-groups})$, we have
 \begin{IEEEeqnarray}{rCl}
  K_0(C_{k+1}^{2n})=\bbz, \quad \qquad  K_0(A_m)=\bbz \oplus \bbz/m\bbz.  \nonumber
 \end{IEEEeqnarray}
Since unitary equivalence gives an isomorphism of the middle $C^*$ algebras and hence an isomorphism of the $K$-groups of middle $C^*$-algebras, 
it follows that the exact sequence \ref{exact} is unitarily equivalent to $\phi_1$ or $\phi_{-1}$. 
This implies that $C_{k+1}^{2n}$
 is isomorphic to $A_1$ or $A_{-1}$. Since $A_1=A_{-1}=C(S_0^{2k+1})$, it follows that $C_{k+1}^{2n}$ is 
 isomorphic to $C(S_0^{2k+1})$. Hence by induction, it follows that  $C(H_q^{2n})$ is 
isomorphic to $C(S_0^{4n-1})$. From Lemma  $3.2$ in \cite{PalSun-2010aa}, it follows that   the $C^*$-algebras $C(S_q^{4n-1})$ are isomorphic 
to $C(S_0^{4n-1})$ for  $q \in (0,1)$. This proves that $C(H_q^{2n})$ is 
isomorphic to $C(S_q^{4n-1})$.
\qed

\brmrk
In case of  $q=0$,  we 
need to be slightly careful to get the defining relations of $C(H_0^{2n})$. In the relation (\ref{c2}), 
we first start with $i=2n$ which gives $z_{2n}z_1=0$. Then we take $i=2n-1$ and so on and get  
 the relation $z_iz_{i^{'}}=0$ for $i < n$. Further in the relation (\ref{c5}), it is easy to check that for $i+j < 2n+1$, $\rho_i +\rho_j >0$.  
Now  by putting $q=0$ in the relations (\ref{c3}), (\ref{c4}) and (\ref{c4}), we get $z_i^*z_j=0$ for $i \neq j$. Other relations 
are obtained by putting $q=0$ in the remaining relations. By looking at the relations, one can see that  the defining relations 
of $C(H_0^{2n})$ are exactly same as those of $C(S_0^{4n-1})$.  
These facts together with  Theorem \ref{isomorphism} prove that for
different values of $q \in [0,1)$, the $C^*$-algebras $C(H_q^{2n})$ are isomorphic.
\ermrk

 \noindent\begin{footnotesize}\textbf{Acknowledgement}:
I would like to thank  Arup Kumar Pal, my supervisor, for his constant support and for his valuable suggestions.
\end{footnotesize}

\noindent{\sc Bipul Saurabh} (\texttt{saurabhbipul2@gmail.com})\\
         {\footnotesize Indian Statistical
Institute, 7, SJSS Marg, New Delhi--110\,016, INDIA}

\end{document}